\documentclass[10pt]{amsart}
\usepackage{palatino, amssymb, amsfonts, latexsym}
\usepackage{amssymb, amsmath, amscd}
\input xy
\xyoption{all}

\newtheorem{theorem}{Theorem}[section]
\newtheorem{lemma}[theorem]{Lemma}
\newtheorem{prop}[theorem]{Proposition}
\newtheorem{cor}[theorem]{Corollary}

\theoremstyle{definition}
\newtheorem{definition}[theorem]{Definition}

\theoremstyle{remark}
\newtheorem{remark}[theorem]{Remark}

\newcommand{\U}{\mathbf{U}_q(\widehat{\mathfrak{sl}}_2)}
\newcommand{\asl}{\widehat{\mathfrak{sl}}_2}
\newcommand{\tp}{\tilde{\psi}}
\newcommand{\UA}{\mathbf U^+_{\mathcal A}}
\newcommand{\Uz}{\mathbf U^+_{\mathcal A}(0)}
\newcommand{\A}{\mathcal A}
\newcommand{\e}{\varepsilon}

\numberwithin{equation}{section}

\begin{document}

\title{The Kronecker quiver and bases of quantum affine $\mathfrak{sl}_2$}

\author{Kevin McGerty}
\address{School of Mathematics, Institute of Advanced Study. }

\date{June 15, 2004}

\begin{abstract}
We compare various bases of the quantum group $\mathbf U (\asl)$ in the context of the Kronecker quiver, and relate them to the Drinfeld presentation.
\end{abstract}

\maketitle

\section{Introduction}\label{quantum}

In this paper we relate for affine $\mathfrak{sl}_2$ the explicit algebraic approximation of the canonical basis in \cite{BCP} to the geometry of the variety of quiver representations. In particular we find that the ``purely imaginary'' elements of the basis in \cite{BCP} are in fact related to the corresponding elements of the canonical basis in a very simple way (even on the algebraic level). We also show how the part of the relations in the Drinfeld presentation may be understood in this context.

We begin by recalling the definition of the quantum group $\U$. This is an algebra over $\mathbb Q(v)$ where $v$ is an indeterminate, generated by elements 

\[
E_i, F_i, K^{\pm 1}_i (i \in \mathbb Z / 2 \mathbb Z) \text{ and } C^{\pm 1/2}, 
\]
subject to the relations
\begin{itemize}
\item $C^{\pm 1/2}$ is central,
\item $C^{1/2}C^{-1/2} = C^{-1/2}C^{1/2} = 1$, $(C^{1/2})^2 = K_0K_1$.
\item $K_iK_i^{-1} = K_i^{-1}K_i= 1$,
\item $K_iE_j= v^{a_{ij}}E_jK_i, \quad K_iF_j = v^{-a_{ij}}F_jK_i$ \quad for all $i,j$,
\item $[E_i, F_j] = \delta_{ij}(K_i - K_i^{-1})/(v-v^{-1})$ \quad for all $i,j$,
\item $\sum_{k=0}^{3}(-1)^{k}E_i^{(k)}E_jE_i^{(3-k)} = 0$ \\
$\sum_{k=0}^{3}(-1)^{k}F_i^{(k)}F_jF_i^{(3-k)} = 0$ \qquad for $i \neq j$.
\end{itemize}
Here $a_{ij}= (-1)^{1+\delta_{ij}}2$, and the divided power in the last relation is
understood in the quantum sense, that is, $x^{(k)} = x^k/[k]!$ with 
\[
[k] = (v^k - v^{-k})/(v-v^{-1})
\]
and $[k]! = [1][2]\ldots[k]$.

The plus part of $\U$ is the subalgebra generated by $\{E_0, E_1\}$, and is denoted $\mathbf{U}^+$. Similarly
we have $\mathbf{U}^-$, the subalgebra generated by $\{F_0, F_1\}$, and $\mathbf{U}^0$ the subalgebra generated by
$\{K_0, K_1\}$. It is known that there is a triangular decomposition
$\U = \mathbf{U}^+.\mathbf{U}^0.\mathbf{U}^-$. Let $\mathcal A = \mathbb Z[v,v^{-1}]$, and let
$\mathbf U^+_{\mathcal A}$ be the $\mathcal A$-algebra generated by $\{ E_i^{(k)} \colon i \in \mathbb Z/ 2 \mathbb Z,
k \in \mathbb N\}$. Then $\mathbf U^+_{\mathcal A}$ is an $\mathcal A$-form of $\mathbf U^+$.
Note also that there is a natural involution $\tau$ of $\U$ which is induced by the map interchanging
the two elements of $\mathbb Z/ 2 \mathbb Z$.

This description of $\U$ is known as the Cartan-Serre presentation. There is an alternative presentation due to Drinfeld \cite{Dr} which is the quantum analogue of the loop realization of an affine Lie algebra. This presentation is important in the study of finite dimensional representations.

\begin{theorem}
Let $A$ be the algebra over $\mathbb Q(v)$ generated by $\{x_r^{\pm}, k^{\pm}, C^{\pm 1/2}, h_s: r \in \mathbb Z,
s \in \mathbb Z \backslash \{0\}\}$
subject to the relations
\begin{itemize}
\item $C^{\pm 1/2}$ is central,
\item $k^+k^- = k^-k^+ = 1, \quad C^{1/2}C^{-1/2} = C^{-1/2}C^{1/2} = 1$
\item $k h_s = h_s k$,
\item $k x_{r}^{\pm} = v^{\pm 2}x_r^{\pm} k$,
\item $[h_r,h_s] = \delta_{r,-s} \frac{1}{r}[ra_{ij}]\frac{C^r - C^{-r}}{v-v^{-1}}$,
\item $[h_s,x_r^{\pm}] = \pm \frac{[2s]}{s} C^{\mp |r|/2}|x_{r+s}^{\pm}$,
\item $x_{r+1}^{\pm}x_s^{\pm} - v^{\pm 2}x_s^{\pm}x_{r+1}^{\pm} = v^{\pm 2}x_r^{\pm}x_{s+1}^{\pm} -
x_{s+1}^{\pm}x_r^{\pm}$,
\item $[x_r^+,x_s^-] = (C^{(r-s)/2}\psi^+_{r+s} - C^{-(r-s)/2}\psi^-_{r+s})/(v-v^{-1})$.
\end{itemize}

where the $\psi_s^{\pm}$ are defined by the relation
\[
    1 + (v-v^{-1})\sum_{s \geq 1} \psi_{\pm s}^{\pm} u^{\pm s} = k^{\pm}\mathrm{exp}\biggl( (v-v^{-1})
    \sum_{r \geq 1} h_{\pm r} u^{\pm r} \biggr).
\]
Then $A$ is isomorphic to $\U$. Moreover the subalgebra generated by $$\{x_r^+, C^{s}k^{-1}x_s^-, h_s \colon
r \geq 0, s > 0 \}$$ is precisely $\mathbf U^+$.
\end{theorem}

The map establishing this isomorphism is nonobvious, and was constructed for a general affine quantum group
using the action of the braid group by Beck in \cite{B}. Later we will see that we can recover ``half'' of the case of $\asl$ by considering the representations of quivers.

\section{The Kronecker quiver}
\label{Kronecker}

In this section we recall the representation theory of the Kronecker quiver, and the quiver approach to quantum groups as studied by Ringel and Lusztig. For more details see \cite{L92}, \cite{L93}. Fix an algebraically closed field $\mathbf k$. The Kronecker quiver $K$ has two vertices, which we shall index by $\mathbb Z / 2\mathbb Z$, and two arrows, both pointing from $0$ to $1$. A representation of $K$ is a pair of vector spaces $V_0, V_1$ and two linear maps $x_1, x_2$ from $V_0$ to $V_1$. Let $\mathcal K$ be the category of finite dimensional representations of $K$ over $\mathbf k$, then $\mathcal K$ is an Abelian category.

There are two simple objects, $S_0$ and $S_1$, where for $S_0$ we have $V_0 = \mathbf k$ and $V_1 = 0$, and for $S_1$ we have $V_0 = 0$ and $V_1 = \mathbf k$. Thus the equivalence class of a representation $V$ in the Grothendieck group of $\mathcal K$ is determined by the vector $\dim(V) = (\dim(V_0),\dim(V_1))$. As a Krull-Schmidt theorem holds in $\mathcal K$, to classify representations it is enough to have a classification of the indecomposables. This was essentially done by Kronecker, as we now describe in modern terminology. There are three different classes of indecomposables: the \textit{preprojective}, \textit{preinjective}
and \textit{regular} modules. The indecomposable preprojective and preinjective  modules are uniquely determined by their dimension vector. For $k \geq 1$ there is a unique indecomposable preprojective module $P_k$ with dimension vector $(k-1, k)$: the maps $x_1, x_2$ are injective and have distinct image. Similarly for each $k \geq 1$ there is a unique indecomposable preinjective $I_k$ of dimension $(k, k-1)$: the maps
$x_1, x_2$ are both surjective, and their kernels are distinct. The regular indecomposables have dimension vector $(n,n)$ for some $n \geq 1$, however there
is a moduli of such representations. Indeed they are parametrized by points of $\mathbb{P}^1(\mathbf k)$ as follows. A representation $((\mathbf k^n, \mathbf k^n), x_1,x_2)$ is regular indecomposable if there exists $[\lambda_0, \mu_0] \in \mathbb P^1(\mathbf k)$ such that $\lambda x_1 + \mu x_2$ is invertible for all $[\lambda, \mu] \neq [\lambda_0, \mu_0]$, and $\lambda_0 x_1 + \mu_0 x_2$ has a one-dimensional kernel. We denote the indecomposable representation of dimension $(n,n)$ corresponding to
$z \in \mathbb P^1$ by $R_{z,n}$. We call a representation regular (respectively preprojective, preinjective) if each of its indecomposable summands is regular (respectively preprojective, preinjective). For a given dimension $\mathbf d = (d_0,d_1)$, and $V=(V_0,V_1)$ a vector space of dimension $\mathbf d$, the isomorphism classes of $\mathbf d$-dimensional representations correspond to the $\mathrm G_{V} = \mathrm{GL}_{V_0} \times \mathrm{GL}_{V_1}$ orbits on
\[
E_{V} =
\mathrm{Hom}(V_0, V_1) \oplus \mathrm{Hom}(V_0, V_1).
\]
Note that $x \in E_V$ is regular precisely if for all but finitely many $[\lambda, \mu] \in \mathbb{P}^1(\mathbf k)$ the map $\lambda x_1 + \mu x_2$ is an isomorphism. Hence the regular representations form an open subset $E^r_{\mathbf d}$ of $E_\mathbf d$ (which is empty unless $\mathbf d= (n,n)$ for some $n \in \mathbb N$).

For later use we record some properties of the trichotomy of representations described above: for $P$ preprojective, $I$ preinjective, and $R$ regular, we have

\begin{equation}
\label{homs}
\begin{split}
\text{Hom}(I,P)=\text{Hom}(I,R)=\text{Hom}(R,P)= 0, \\
\text{Ext}^1(P,R)= \text{Ext}^1(R,I) = \text{Ext}^1(P,I) = 0.
\end{split}
\end{equation}
In fact it is also possible to show that 
\begin{equation}
\label{exts}
\text{Ext}^1(I_j,I_k) = \text{Ext}^1(P_k,P_j) = 0 \qquad (j \geq k).
\end{equation}

We now describe Lusztig's construction of the canonical basis. Fix a prime $l$ coprime to the characteristic of $\mathbf k$. For a variety $X$ over $\mathbf k$, we write $\mathcal D(X)$ for the bounded derived category of complexes of $\overline{\mathbb Q}_l$ vector spaces on $X$. We will use the notation of \cite{BBD} for perverse sheaves etc. and will use the term semisimple complex for a complex which is isomorphic to a direct sum
of perverse sheaves with shifts.
For $\mathbf d$ in $\mathbb N^2$, let $S_{\mathbf d}$ be the set of all pairs $(\mathbf i, \mathbf a)$ where $\mathbf i = (i_1, i_2, \ldots, i_m)$ is a sequence of elements in $\mathbb Z/ 2 \mathbb Z$, and $\mathbf a = (a_1,a_2,\ldots,a_m)$ is a sequence of integers, so that $\sum a_k i_k = d_0 0 + d_1 1 $ (in the group algebra of $\mathbb Z/2 \mathbb Z$,
which we identify with the Grothendieck group of $\mathcal K$).
If $\mathrm V_{\mathbf d}$ is a space of dimension $\mathbf d$, a flag of type $(\mathbf i, \mathbf a)$ is a sequence of subspaces $\mathfrak f = (V_\mathbf{d} = V^0 \supset V^1 \supset \cdots \supset V^m = 0)$ such that $\dim(V^k/V^{k+1}) = a_i i_k$. Let $\mathfrak F_{\mathbf i, \mathbf a}$ be the space of all flags of type $(\mathbf i, \mathbf a)$, let $\tilde{\mathcal F}_{\mathbf i, \mathbf a}$ be the variety of pairs $\{(x, \mathfrak f) \in E_{\mathbf d} \times \mathfrak F_{\mathbf i, \mathbf a} \colon x$ preserves $\mathfrak f\}$, and let $\pi_{\mathbf i, \mathbf a}:\tilde{\mathcal F}_{\mathbf i, \mathbf a} \to E_{\mathbf d}$ be the obvious $G_\mathbf d$-equivariant map.
By the decomposition theorem, the complex $L_{\mathbf i, \mathbf a} = (\pi_{\mathbf i, \mathbf a})_!(1)$ is a direct sum of simple perverse sheaves with shifts, where $1$ is the constant complex on $\tilde{\mathcal F}_{\mathbf i, \mathbf a}$. Let $\mathcal P_{\mathbf d}$ be the (finite) set of simple perverse sheaves that occur with some shift in some complex $L_{\mathbf i, \mathbf a}$, and $\mathcal Q_{\mathbf d}$ the subcategory of $\mathcal D(E_\mathbf{d})$ consisting of complexes isomorphic to a direct sum of $P[l]$ for various $P \in \mathcal P_{\mathbf d}$, and various $l \in \mathbb Z$. (It is clear that if $V'$ is another vector space of dimension $\mathbf d$ then the simple perverse sheaves obtained by this construction on $E_{V'}$ are canonically isomorphic to those on $E_V$, hence the notation $\mathcal P_\mathbf d, \mathcal Q_\mathbf d$ is justified.)

Now take $\mathbf b$ and $\mathbf c \in \mathbb N^2$, such that $\mathbf b + \mathbf c = \mathbf d$,
and let $V_{\mathbf b}, V_{\mathbf c}$ be vector spaces of dimensions $\mathbf b, \mathbf c$
respectively. Consider the diagram:

\xymatrix{
& & & E_{V_{\mathbf b}} \times E_{V_{\mathbf c}}& \ar[l]_{p_1} E' \ar[r]^{p_2} & E'' \ar[r]^{p_3} &
E_{V_{\mathbf d}}.}
\noindent
Here $E'$ is the variety of triples $(x,\phi, \psi)$ where $x \in E_{\mathbf d}$,

\xymatrix{& & & 0 \ar[r] & V_{\mathbf b} \ar[r]^{\phi} & V_{\mathbf d} \ar[r]^{\psi} & V_{\mathbf c} \ar[r] & 0}

\noindent
is an exact sequence of $\mathbb Z/ 2\mathbb Z$-graded vector spaces, and the image of $\phi$ is $x$-stable.
The variety $E''$ consists of pairs $(x, W)$ where $W$ is an $x$-stable subspace of $V_{\mathbf d}$ with $\dim(W) = \mathbf b$. The map $p_1$ is given by $(x, \phi, \psi) \mapsto (\phi^{-1} x \phi, \psi \bar{x} \psi^{-1})$ where $\bar{x}$ is the map induced by $x$ on $V_{\mathbf d}/ \rm{im}(\phi)$, while $p_2$ is given by $(x, \phi, \psi)
\mapsto (x, \rm{im}(\phi))$ and $p_3$ is given by $(x, W) \mapsto x$. One notes that $p_1$ is smooth, while $p_2$ is a principal $G_{V_\mathbf b}\times G_{V_\mathbf c}$ bundle, and $p_3$ is proper.

Given complexes $L$ and $L'$ in $\mathcal Q_{\mathbf b}$ and $\mathcal Q_{\mathbf c}$ respectively we consider the exterior product $L\boxtimes L'$. Now the map $p_1$ is smooth and $G_{\mathbf c_1}\times G_{\mathbf c_2}$-equivariant, so $p_1^*(L\boxtimes L')$ is a direct sum of simple perverse sheaves with shifts, equivariant under the action of
$G_{\mathbf b}\times G_{\mathbf c}$. As $p_2$ is a principal bundle, there is a semisimple complex $A$ on $E''$ such that $p_2^*(A) \cong p_1^*(L\boxtimes L')[m](m)$ (where $(\cdot)$ denotes the Tate twist). Here $m$ is the difference of the fibre dimensions of $p_1,p_2$. More explicitly,
\[
m = m(\mathbf b,\mathbf c) = 2b_0c_1 + b_0c_0 + b_1c_1.
\]
Set $L*L' = (p_3)_!(A)$. It is not hard to see that this maps $\mathcal Q_{\mathbf c_1}\times \mathcal Q_{\mathbf c_2}$ to $\mathcal Q_{\mathbf d}$. Let $\mathbf K_{\mathbf d}$ be the Grothendieck group of $\mathcal Q_{\mathbf d}$, it is a $\mathcal A$ module via
$v(L) = L[1]$ and $v^{-1}(L) = L[-1]$. Then $\mathbf K = \bigoplus_{\mathbf d \in \mathbb N^2} \mathbf K_{\mathbf d}$ becomes an associative $\mathcal A$ algebra under $*$ (in the Hall algebra context of Ringel this corresponds to the twisted Hall algebra of the Kronecker quiver).

If $W = (W_0,W_1)$ is a fixed subspace of $V$ of dimension $\mathbf c$, then consider the set $F$ consisting of those $x \in E_V$ which preserve $W$. There is a natural map from $F$ to $E_T \times E_W$ where $T= V/W$, and hence we have a diagram:

\xymatrix{& & & & E_T \times E_W & \ar[l]_{\kappa} F \ar[r]^\iota & E_V.}
\noindent
Then set $\Delta_{\mathbf b, \mathbf c}(A) = \kappa_!(\iota^*(A))[n](n)$,
where $$n= n(\mathbf b, \mathbf c) = 2b_0c_1 - b_0c_0 - b_1c_1.$$
It is shown in \cite{L93} that this maps $\mathcal Q_\mathbf d$ to $\mathcal Q_{\mathbf b,\mathbf c}$, where the latter is defined in the obvious way with respect to two distinct copies of the Kronecker quiver. Summing over all $(\mathbf b, \mathbf c)$ such that $\mathbf b + \mathbf c = \mathbf d$ we obtain a map $\Delta: \mathbf K \to \mathbf K\otimes \mathbf K$, which gives $\mathbf K$ a ``twisted'' coalgebra structure.

Let $\theta_i$ be the constant complex on $E_{i}$ for $i \in \mathbb Z/ 2\mathbb Z$ (note that $E_i$ is a point). We have the following theorem due to Lusztig:

\begin{theorem}
\label{quiver}
There is a unique $\mathcal A$-algebra isomorphism $\chi : \mathbf K \to \mathbf U^+_{\mathcal A}$ mapping $\theta_i$ to $E_i$.
\end{theorem}

See \cite{L93} for a proof of this for a general quantum group. Via the map $\chi$, the simple perverse sheaves
 in $\mathcal P_{\mathbf d}$ provide $\mathbf U^+_\mathcal A$ with a natural basis -- the canonical basis. In \cite{L92} Lusztig described the elements of the canonical basis explicitly for affine quantum groups, by giving the support of the irreducible perverse sheaves and the local systems they restrict to on an open dense locus. We will give this description in the case of $\asl$. Let $\mathcal S^1_\mathbf d$ be the
set of pairs $((r_i),(s_i))$ where $(r_i)$ and $(s_i)$ are sequences of nonnegative integers such that
\[
\sum_{i \geq 1} r_i\dim(P_i) + \sum_{i \geq 1} s_i\dim(I_i) + p(1,1) = \mathbf d,
\]
for some integer $p \geq 0$.

For $\sigma \in \mathcal S^1_{\mathbf d}$ let $X(\sigma)$ denote the subset of $E_\mathbf d$ which correspond to representations of the Kronecker quiver which are isomorphic to
\[
r_1P_1 \oplus r_2P_2 \oplus \ldots \bigoplus s_1I_1 \oplus s_2I_2 \oplus \ldots \bigoplus R_{z_1,1}\oplus R_{z_2,1}
\oplus \ldots R_{z_p,1},
\]
for some distinct $z_i \in \mathbb P^1$. Let
\[
\begin{split}
\tilde{X}(\sigma) = \{(x,z_1,z_2,\ldots,z_p): x \in X(\sigma), x \text{ contains} \\ \text{a submodule isomorphic to }
 R_{z_i,1}\}
\end{split}
\]
There is an obvious map $\tilde{X}(\sigma) \to X(\sigma)$ which is a principal $S_p$ covering ($S_p$ the symmetric group on $p$ letters).
Thus for every partition $\lambda$ of $p$ we have a local system $\mathcal L_{\lambda}$ on $X(\sigma)$, and hence
an intersection cohomology sheaf on $\overline{X(\sigma)}$. Let $\mathcal S_{\mathbf d}$ be the set of triples
$((r_i),(s_i),\lambda)$ such that
\[
\sum_{i \geq 1} r_i\dim(P_i) + \sum_{i \geq 1} s_i\dim(I_i) + \sum_{i \geq 1}\lambda_i(1,1) = \mathbf d.
\]

\begin{theorem}
\label{description}
\cite{L92} $\mathcal P_d$ is parameterized by the set $\mathcal S_{\mathbf d}$. Given $(\sigma, \lambda)
\in \mathcal S_{\mathbf d}$ the simple perverse sheaf corresponding to it is the intersection
cohomology extension of the local system $\mathcal L_{\lambda}$ on $\overline{X(\sigma)}$
\end{theorem}

Though we will not prove this theorem here, note that the proof of Theorem \ref{kostka} will at least show that the elements corresponding to $(0,0,\lambda)$, where $\lambda$ is a partition of  $n$, lie in $\mathcal P_{(n,n)}$. Moreover, the fact that the theorem gives the right number of basis elements is immediate from the PBW theorem and the generating function identity:

\[
 \prod_{i \geq 1}(1-x^i)^{-1} = \sum_{i \geq 0} p(i)x^i,
\]
where $p(i)$ is the number of partitions of $i$ (and $p(0)$ is understood to be $1$).

\section{Root Vectors and the Hall algebra}

In this section we will work in the context of the Hall algebra over a finite field following the account given in \cite{L98} whose notation we will also follow. \cite[\S 5]{L98} describes the elements of the canonical basis in this context, which we shall need later. Thus let $\mathbf k$ be an algebraic closure of the field $\mathbb F_p$ equipped with the action of a Frobenius $F$ such that $\mathbf k^F = \mathbb F_q$ the finite field with $q$ elements. For a variety $X$ over $\mathbf k$ let $X^F$ denote the fixed points of the Frobenius on $X$. The analogue of $\mathbf K$ is an algebra $\mathcal F$ of $\overline{\mathbb Q}_l$-valued functions on the various $E_V^F$. Note that $v$ in the definition of $\mathbf U^+_{\mathcal A}$ can be specialized to any non-zero $\varepsilon \in \overline{\mathbb Q}_l$. Denote the resulting algebra by $\mathbf U^+_\varepsilon$, and let $\mathcal A_\varepsilon$ be the subring $\mathbb Z[\varepsilon, \varepsilon^{-1}]$ in $\overline{\mathbb Q}_l$. We will use the notation $[n]_\e$ for $[n]_{|v= \e}$ etc. If we fix $\varepsilon$ to be a square root of $q$, then \cite{L98} shows that there is an isomorphism of $\mathcal A_\e$-algebras, $\chi_\varepsilon: \mathcal F \to \mathbf U^+_\varepsilon$.

In \cite{BCP} an integral PBW basis for the positive part of an affine quantum group is constructed. We wish to show that the root vectors from which this basis is constructed have a very natural description in the context of the Hall algebra.

We briefly describe their construction. The root vectors corresponding to real roots are constructed using the braid group as in \cite{L91}. Recall that the braid group $B$ of $\asl$ is freely generated by two elements $T_0, T_1$. It acts on $\U$ in the standard way (see \cite{L93} for a detailed discussion of this), and restricts to an action on the integral form of the quantum group, as is clear from the explicit formulas defining the action. In the case of $\asl$ we have $T_0\tau = \tau T_1$, where $\tau$ is the involution defined in section \ref{quantum}.

Let $T_{\varpi}$ denote this automorphism, and set
\[
E_{0,n} = T_{\varpi}^{n}(E_0) = T_0T_1\ldots T_{n-1}(E_n),
\]
and
\[
E_{1,n}= T_{\varpi}^{-n}(E_1) = T_1^{-1}T_0^{-1}\ldots T_{n}^{-1}(E_{n+1})
\]
where $n-1, n$ and $n+1$ should be understood modulo $2$.
Next we define elements $\tilde{\psi}_k$ by setting
\[
    \tp_k = E_{0,k-1}E_1 - v^{-2}E_1E_{0,k-1}.
\]
Finally we recursively define elements $\tilde{P}_k$ by setting $\tilde P_0 = 1$ and setting
\[
\tilde{P}_k = \frac{1}{[k]}\sum_{r=1}^k v^{r-k} \tp_r \tilde{P}_{k-r},
\]
The PBW basis is then constructed by taking products of these elements using a suitable ordering.

Note that while the elements $\tp_k$ clearly lie in the integral form $\UA$, it is not obvious from the above definition that the same holds for the elements $\tilde{P}_k$. Nevertheless this is shown in \cite{BCP} using the results of \cite{CP}.

The aim of the rest of this section is to describe the specializations $E_{i,n}, \tilde{P}_k, \tp_k \in \mathbf U_\e^+$  as elements of $\mathcal F$ via the isomorphism of Theorem \ref{quiver}. The simplest of these to describe are the elements $E_{i,n}$, since they are defined using the action of the braid group, and a description of this action on $\mathcal F$ is known. In the finite type case it was described in \cite{L91}, and in full generality in \cite{L98}, we review it here in the context we need.

We need to introduce another orientation of the Kronecker quiver -- from now on we will write $K_{+}$ for the quiver with two arrows pointing from $0$ to $1$, and $K_{-}$ for the quiver with two arrows pointing from $1$ to $0$. We will write $E_{V, \pm}$ for the space of representations of $K_{\pm}$ with on $V$, and use similar modifications of the notation of Section \ref{Kronecker} when we need to specify the orientation we are using. Let $^1E_{V,+}$ be the open subset of $E_{V,+}$ consisting of those $x \in E_{V, +}$ such that 
\[
x_1(V_0) + x_2(V_0) = V_1,
\]
and let $E^1_{V,-}$ be the open subset of $E_{V,-}$ consisting of those $x \in E_{V,-}$ such that 
\[
\mathrm{ker}(x_1) \cap \mathrm{ker}(x_2) = 0.
\]
If $\dim(V)=\mathbf d$ let $\mathbf d' = (d_0,2d_0-d_1)$, and let $V'$ be a vector space of dimension $\mathbf d'$, with $V_0 = V'_0$. Define $Z$ to be the set of pairs $(x^{-}, x^{+}) \in$ $E^1_{V, -} \times$ $^1E_{V',+}$ such that we have an exact sequence

\xymatrix{& & & 0 \ar[r] & V_1 \ar[r]^{(x^-_{1},x^-_{2})} &
 V'_0 \oplus  V'_0 \ar[r]^{x^+_1+ x^+_2} & V'_1 \ar[r] & 0.}
\noindent Hence we have a diagram

\xymatrix{& & & & E^1_{V, -} & \ar[l]_\alpha Z \ar[r]^{\beta} & ^1E_{V',+}}

\noindent
where $\alpha$ is a principal $GL_{V'_1}$ bundle and $\beta$ is a principal $GL_{V_1}$ bundle. 

Let $j^-_1 \colon E^1_{V, -} \to E_{V, -}$ and let $j^+_1: $ $^1E_{V',+} \to E_{V',+}$. Let $\mathcal F^1_{V,-}$ denote the functions in $\mathcal F_V$ supported on $E^{1, F}_{V,-}$ and $^1\mathcal F^F_{V',+}$ the functions in $\mathcal F_{V'}$ supported on $^1E^F_{V, +}$.

\begin{definition}
Let $\sigma_1: \mathcal F^1_{V,-} \to {^1\mathcal F}_{V',+}$ be defined as follows. Given $f \in
\mathcal F^1_{V,-}$, there is a unique function $f' \in {^1\mathcal F}_{V', +}$ such that
$\alpha^*(f) = \beta^*(f')$. Set

\[
\sigma_1(f) = \e^{\dim(GL_{V_1'}) - \dim(GL_{V_1})} f'.
\]
\end{definition}

In a completely analogous way we may define $^0\mathcal {F}_{W,-}$, and $\mathcal F^0_{W,+}$, and a
map $\sigma_0: \mathcal F^0_{W,+} \to {^0\mathcal F}_{W',-}$ where if $\dim(W) = \mathbf c = (c_0, c_1)$
then $\dim W' = \mathbf c' = (2c_1-c_0, c_1)$

\begin{remark}
We remark that in the context of perverse sheaves since
the maps $j^{\pm}_1$ are open embeddings, there is a bijective correspondence (essentially the intersection
cohomology extension) between simple perverse sheaves on, say, $E_{V,-}$ whose supports have dense intersection
with $E^1_{V,-}$, and simple perverse sheaves on $E^1_{V,-}$. Clearly this extends to an
injection from the class of semisimple perverse sheaves on $E^1_{V,-}$ to semisimple perverse sheaves
on $E_{V,-}$.
Moreover, as the maps $\alpha, \beta$ are principal bundles, the choice of shift ensures that the analogue of
$\sigma_1$ on the level of perverse sheaves gives a bijection between simple perverse sheaves on $E^1_{V,-}$
and simple perverse sheaves on $^1E_{V', +}$. This gives a geometric incarnation of \cite[Theorem 1.2]{L96},
certainly known to Lusztig, but perhaps not written in the literature.
\end{remark}

The choice of orientation of the quiver $K$ is known not to effect the validity of Theorem \ref{quiver},
and so $\sigma_1$ induces a map $\varrho_1$ from $\mathbf U^+_\varepsilon$ to itself. Indeed in our case it is not
even necessary to appeal to this result -- if write $\tau$ for the obvious map
from $E_\mathbf d$ to $E_{\mathbf d^t}$ where $\mathbf d^t = (d_1, d_0)$, then $\tau$ induces
an isomorphism from $\mathcal F_+$ to $\mathcal F_-$, and if $\chi^{\pm}_\e$ denote the isomorphisms from
$\mathcal F_{\pm}$ to $\mathbf U^+_\e$ then $\chi^+_\e \tau = \tau \chi^-_\e$ where $\tau$ on the left side of the equation
is the map just defined, and $\tau$ on the right side is the map defined in Section \ref{quantum}.

\begin{prop}\cite[\S8,\S9]{L98}
The map $\varrho_1$ coincides with $T_1$ on the subalgebra $\{ x \in \mathbf U^+_{\varepsilon}: T_1(x) \in
\mathbf U^+_{\varepsilon}\}$
\end{prop}

Now we may easily identify the elements $E_{i,n}$. It is straight forward to check the following.

\begin{prop}
\label{orbits}
We have
\begin{enumerate}
\item Let $\mathbf d = (k, k-1)$, $V$ a space of dimension $\mathbf d$ and let $\mathcal O_\mathbf d$
be the open dense orbit of $G_V$ on $E_V$ and let $\gamma_k =
\e^{-\dim(E_V)}1_{\mathcal O^F_{\mathbf d}} \in \mathcal F^0_{V, +}$, then $\chi(\gamma_k) = E_{0,k}$.

\item Let $\mathbf e = (k-1, k)$ and let $W$ be a space of dimension $\mathbf e$. Then let
$\mathcal O_\mathbf e$ be the open dense orbit of $G_W$ on $E_W$.
Let $\mu_k = \e^{-\dim (E_W)} 1_{\mathcal O^F_\mathbf e}$,
then $\chi(\mu_k) = E_{1,k}$.
\end{enumerate}
\end{prop}

\begin{remark}
The dense orbits in $(1)$ and $(2)$ of the previous proposition correspond to the indecomposables $I_k$,
respectively $P_k$, described in Section \ref{quiver} -- one can use the description of these representations
to check the assertion that the orbits are indeed dense. Note also that this shows that the functions
$1_{\mathcal O^F_{\mathbf d}}$ and $1_{\mathcal O^F_\mathbf e}$ are in $\mathcal F$, which is not
immediately obvious.
\end{remark}

Finally we wish to identify the elements $\tilde P_k$. Let $\phi_k$ denote the elements of $\mathcal F$ corresponding to $\tp_k$. Recall that $E_V^r$ is the set of regular elements, and that $E_V^r$ is open and nonempty if and only if $\dim(V) = (k,k)$ for some $k \geq 1$. In that case let $j^r: E_V ^r \to E_V$, and let $\rho_d = \e^{-\dim(E_V)} 1_{E_V^r}$ (with $\rho_0 = 1$).

\begin{lemma}
\label{reg}
The functions $\rho_d$ are in $\mathcal F$.
\end{lemma}

\begin{proof}
We show this by induction on $d$. For $d=1$ this immediate from the direct computation
\[
\rho_1 = \theta_0\theta_1 - q\theta_1\theta_0,
\]
where $\theta_i$ is the indicator function of $E_W$, $\dim(W)= i \in \mathbb Z/2\mathbb Z$.
For $d \geq 1$ we use induction. Suppose $\rho_k \in \mathcal F$ for all $k \leq d$.
It is easy to check that $\theta_0^{(d)}\theta_1^{(d)} = 1_{E_W}$. Thus we need to show
that the function $1_{E_V \backslash E_V^r}$ lies in $\mathcal F$.
To do this note first that if $x \in E_V \backslash E_V^r$ then the representation $(V,x)$
has a preprojective and/or preinjective component. Indeed it follow from the representation theory of $K$ that $(V,x)$ can be written uniquely as a sum $P \oplus R \oplus I$ where $P,R,I$ are preprojective, preinjective, and regular respectively, and if $x$ is not regular, then at least one of $P$ and $I$ is nonzero.
The functions $\gamma_i$ are in $\mathcal F$ and so it follows that if $W$ is a
$\dim(I)$-dimensional vector space, and $\mathcal O_I$ is the orbit of $G_W$ corresponding
to $I$, then $\gamma_I = \e^{-\dim(O_P)/2}1_{\mathcal O_P}$ lies in $\mathcal F$. Indeed if $I \cong I_{i_1} \oplus I_{i_2} \oplus \cdots \oplus I_{i_r}$ say, with $i_1 < i_2 < \cdots i_r$, then one can check using \eqref{exts} that $\gamma_I = \e^{\alpha} \gamma_{i_r} \gamma _{i_{r-1}} \ldots \gamma_{i_1}$ where $\alpha \in \mathbb Z$. Similarly if we let $\mu_P = \e^{-\dim(\mathcal O_P)}1_{\mathcal O_P}$ where $\mathcal O_P$ is the orbit corresponding to $P$, and $P \cong P_{j_1} \oplus P_{j_2} \oplus P_{j_3} \oplus \cdots \oplus P_{j_s}$ then we have $\mu_P = \e^{\beta} \mu_{j_1} \mu_{j_2} \ldots \mu_{j_s} \in \mathcal F$ for $\beta \in \mathbb Z$.

Finally, it is then clear using \eqref{homs} that there exist constants $c_{P,I}$ such that 
\[
\e^{-\dim(E_V)}1_{E_V \backslash E_V^r} = \sum_{k < d} c_{P,I}\mu_P r_k \gamma_I,
\]
where the sum is over all triples $P,k,I$ where $P$ is preinjective, $I$ preprojective and $k$ is strictly less than $d$.
\end{proof}

It follows from the representation theory of $K$ that the subalgebra $\mathcal F^r$ consisting of functions in $\mathcal F$ which are supported on $E_V^r$ is commutative. Let us define, for $f \in \mathcal F_V$ the function $r(f)$ to be the resriction of $f$ to $E_V^r$, that is $r(f) = (j^r)^*(f)$. Note that it is not \textit{a priori} clear that this restriction is necessarily in $\mathcal F$ itself, though we will later see that this is indeed the case. However it is easy to see that the functions $\phi_k$ are in $\mathcal F^r$, and that $\phi_k = r(\gamma_{k-1} \theta_1)$.

\begin{theorem}
\label{relation}
In $\mathcal F$ we have
\[
[n]\rho_n = \sum_{i=1}^{n} \e^{i-n}\phi_i \rho_{n-i},
\].
\end{theorem}

\begin{proof}
This is essentially contained in \cite[Theorem 4.1]{Z}, but we give a different, more conceptual,
proof here. Note that $\sum \e^{i-n}\phi_i \rho_{n-i}$ is regular ($\mathcal F^r$ is a
commutative subalgebra of $\mathcal F$), and so we need only compute its value
at regular elements. It is more convenient to compute the product in the opposite order:
\begin{equation}
\label{terms}
\sum \e^{i-n}\rho_{n-i}\phi_i= \sum \e^{i-n}\rho_{n-i}\gamma_{i-1}\theta_1 -
\sum \e^{i-n-2}\rho_{n-i}\theta_1\gamma_{i-1}
\end{equation}
It follows immediately from \eqref{homs} that the second term on the righthand side restricts to
zero on the regular set. Thus it remains to compute the first term at a point $x \in E_V^r$. A simple calculation with the cocycle $m(\cdot,\cdot)$ show that this value is just $\e^{-2n^2 -n +1}$ times the number of filtrations of $V$ of the form $(L < W < V)$ where $L \cong S_1$ and $W= (W_0,W_1)$ has $\dim(W) = (k,k)$ ($1 \leq k \leq n$), such that the quotient $W/L$ is isomorphic to $I_{k-1}$ and the quotient $V/W$ is regular. Now it is clear that the submodule $W$ must be regular, and hence the quotient $V/W$ is automatically regular, as the regular modules for an Abelian category. Thus the value of the first term on the right in Equation \eqref{terms} at $x \in E_V^r$ is 

\begin{equation}
\label{counting}
\begin{split}
\e^{-2n^2 -n + 1} |\{(W,L): L < W < V, L \cong S_1, W \text{ regular}, \\ W/L \text{ indecomposable preinjective}\}|,
\end{split}
\end{equation}

For any line $L \subset V_1$, there is a unique minimal regular submodule $R$ of $(V,x)$ which contains $L$. Let $(j,j)$ be the dimension of $R$. Then $(V/L,x)$ contains $R/L$ as a submodule. The point is to observe that $R/L$ is isomorphic to $I_{j-1}$, i.e. that $x_1^{-1}(L)$ and $x_2^{-1}(L)$ are distinct lines in $R_0$. 

The minimality of $R$ ensures that these subspaces are disjoint, since a line in their intersection would, in $R$, map onto $L$ (as otherwise $R$ could not be regular) and thus provide a smaller regular submodule containing $L$. So it only remains to show that they are one-dimensional. 
Consider the subspace $L + \rm{im}(x_1)$ of $R_1$. It is a proper subspace if and only if $x_1^{-1}(L)$ is not a line. Indeed if it is proper, since $\dim(R_0) = \dim(R_1)$ the kernel of $x_1$ must be either at least one-dimensional, with the image containing $L$, or of dimension at least two, with $L$ not necessarily contained in the image of $x_1$. In either case, the preimage of $L$ is at least two-dimensional. The converse is similar. Picking a subspace of $R_0$ of dimension $\dim(L + \mathrm{im}(x_1))$ which is mapped into $L + \mathrm{im}(x_1)$ by $x_2$ we obtain a proper regular submodule of $R$ which contains $L$, contradicting the minimality of $R$. (The fact that such a subspace exists follows from the same considerations as above). Similarly we see that $x_2^{-1}(L)$ is also a line.

Now for any pair $(L,W)$ as in \eqref{counting} the minimal regular module $R$ considered above certainly lies in $W$ and so $W/L$ contains $R/L$ as a submodule. But an indecomposable preinjective does not contain any indecomposable preinjective of smaller dimension as a submodule, and so $W/L$ is an indecomposable preinjective precisely when $W=R$. Thus we see that the pair $(W,L)$ is determined by $L$, and every line $L$ occurs in some pair. It follows that the number of pairs we wish to count is precisely $(q^n-1)/(q-1)$ or $\e^{n-1}[n]_\e$. Combining this with \eqref{counting} the result follows.
\end{proof}

We now list some consequences of this last result. First note that it can be neatly expressed in
terms of generating functions: let $\wp(u) = \sum_{i \geq 0} \rho_k u^k$ and $\Phi(u) =
\sum_{i \geq 1} \phi_k u^k$. It is straightforward to see that Theorem \ref{relation}
can be rewritten in the form
\begin{equation}
\label{pseries}
\wp(\e^{-1}u)/\wp(\e u)= 1 +(\e-\e^{-1})\Phi(u).
\end{equation}
This demonstrates that the $\phi_k$ are in the $\mathcal A_\e$-algebra generated by the $\rho_k$, since the power series $\wp(u)$ is a unit in the ring of formal power series over $\mathcal F^r$. This is essentially the content of \cite[\S3]{Z}.

\begin{cor}
\label{properties}
We have:
\begin{enumerate}
\item $\chi_v(\rho_k) = \tilde{P}_k$, hence the $\tilde{P}_k$ lie in $\mathbf U^+_\e$.
\item The elements $\rho_k$ are fixed by $T_\varpi$,
\item the elements $\phi_k$ are fixed by $T_\varpi$,
\item for all $k,l \geq 0$ we have
\[
\gamma_k \mu_l - q^{-1}\mu_l \gamma_k = \psi_{k+l-1}.
\]

\end{enumerate}
\end{cor}

\begin{proof}
Since by definition $\chi_\e(\phi_k) = \tilde{\psi}_k$, and $\rho_0 =1$, the previous theorem
shows that $\chi_\e(\rho_k) = \tilde{P}_k$.
It is easy to check directly from the definition that $\rho_k$ is preserved by $T_\varpi$, for all $k \geq 1$.
Then as \eqref{pseries} shows that the $\psi_k$ are in the subalgebra generated by the $\rho_k$, it follows that they too are preserved by $T_\varpi$. The last equation is an immediate consequence of $(2)$ and Proposition \ref{orbits}.
\end{proof}

\begin{remark}
Lemma \ref{reg}, Theorem \ref{relation}, and $(4)$ of it's corollary are all proved in \cite{Z96}, \cite{Z}, and the proofs here of course have some similarities to those papers, however they are somewhat more conceptual, and take advantage of the braid group action which is not used in \cite{Z96}, \cite{Z}. The idea of taking the ``regular part'' of a module is also studied in those papers and corresponds to our map $r$ (indeed we have used the same notation). Our main point here is that the algebraic constructions of \cite{BCP} have a very simple interpretation in the quiver context.
\end{remark}

\begin{remark}
The elements $\tilde{P}_k$ act on a the ``$l$-highest weight''  vector of a
finite dimensional representation of $\U$ via the Drinfeld polynomials, thus it is
intriguing that they have such a simple relation to the canonical basis in this case.
Interestingly they also have a nice description in the context of the Hall algebra of
coherent sheaves on $\mathbb P^1(\mathbb F_q)$ which has been studied by Kapranov and
Baumann-Kassel. This produces a different subalgebra of the affine quantum group which
however still contains the $\tilde{P}_k$ -- see \cite{BK}.
\end{remark}

\begin{remark}
\label{generic}
Note that the results of this section combine to give an explicit description in $\mathcal F$ of the elements used in \cite{BCP} to generate a PBW basis of $\UA$. Since $\mathcal F$ gives only a specialization of $\UA$, we wish to briefly recall how the generic case may also be recovered. Let $\mathcal S$ be an infinite set of finite fields, say all finite fields of characteristic $p$, and consider the various algebras $\mathcal F^s$, $s \in \mathcal S$, with $\e_s = (\sqrt{p})^e$ where $s = \mathbb F_{p^e}$. Let $\mathcal O = \prod_{s \in \mathcal S} \mathcal F^s$, an algebra over $\overline{\mathbb Q}_l$. Then $\mathcal O$ is an $\mathcal A$ algebra
via the map which sends $v$ to the element whose $s$-th component is $\e_s$. We denote by  $\theta_i$ the element whose components are $\theta_i \in \mathcal F^s$. Let $\mathcal F_{\mathcal S}$ be the $\mathcal A$-subalgebra of $\mathcal O$ generated by the elements $\theta_i$. It is shown in \cite{L98} (see also \cite{Gr}) that the obvious map from $\mathbf U^+$ to $\mathcal F_{\mathcal S}$ is an isomorphism. All of the results of this section apply to $\mathcal F_{\mathcal S}$, and hence to $\UA$. Thus for example we recover the result that
the elements $\tilde{P}_k$ are integral.
\end{remark}

\section{The approximate canonical basis}
Next we examine the construction used in \cite{BCP} used to produce an approximation to the canonical basis (there called a ``crystal basis'') from the PBW basis. They define, for a partition $\lambda = (\lambda_1, \lambda_2, \ldots, \lambda_p)$, following \cite{M}
\[
s_\lambda = \det (\tilde{P}_{\lambda_i-i+j})_{1\leq i,j \leq t},
\]
where $t$ is at least the length of the partition $\lambda$. Let $\Uz$ be the subalgebra of $\UA$ generated by the $\tilde P_k$. It is proved in \cite{BCP} that the elements $\tilde P_k$ are algebraically independent and so we may identify $\Uz$ with the algebra of symmetric functions by mapping the elements $\tilde{P}_k$ to the complete symmetric functions, then this maps the $s_{\lambda}$ to the Schur functions. We set $\tilde{P}_\lambda = \tilde{P}_{\lambda_1}\tilde{P}_{\lambda_2}\ldots \tilde{P}_{\lambda_p}$, and $\rho_\lambda = \rho_{\lambda_1}\rho_{\lambda_2}\ldots\rho_{\lambda_p}$.

The next result shows that the purely imaginary elements of the ``crystal basis'' of \cite{BCP} have a very simple description in terms of the canonical basis. For a partition $\lambda$ of $n$, let $b_\lambda$ be the canonical basis element in $\mathcal F$ corresponding to the element $(0, \lambda) \in \mathcal S_{(n,n)}$.

\begin{prop}
\label{kostka}
Let $\mathbf d = (n,n)$, and let $\lambda$ be a partition of $n$, then we have
\[
\rho_{\lambda} = \bigoplus_{\mu \geq \lambda} K_{\mu \lambda} r(b_\mu)
\]
where $K_{\mu \lambda}$ is the Kostka number associated to the partitions $\lambda$ and $\mu$.
\end{prop}

\begin{proof}
Consider the variety $Y  = \mathfrak F_{\mathbf i,\mathbf a}$ where
\[ \mathbf i = (0,1,0,1,\ldots,0,1), \qquad
\mathbf a = (\lambda_1,\lambda_1, \ldots, \lambda_p,\lambda_p).
\]
Let $U$ be the dense open set of $x \in E_{(n,n)}$ whose spectrum consists of $n$ distinct points in $\mathbb{P}^1$. It is straightforward to see that over $U$ the map $\pi = \pi_{\mathbf i,\mathbf a}$ is a quotient of a principal $S_n$ covering corresponding to the subgroup $S_{\lambda_1} \times S_{\lambda_2} \times \cdots \times S_{\lambda_p}$. Standard properties of coverings and the representation theory of the symmetric group show that
\[
\pi_!(1)\mid_U = \bigoplus_{\mu \geq \lambda} K_{\mu, \lambda} \mathcal L_\mu,
\]
where $\mathcal L_{\mu}$ is the local system on $U$ associated the the irreducible representation of $S_n$ indexed by $\mu$, and $K_{\mu, \lambda}$ is the Kostka number associated to the partitions (see \cite[\S 1.7]{M} for example). Since $\pi$ is proper the decomposition theorem tells us that the intersection cohomology complexes  $IC(\mathcal L_\lambda)$ corresponding to these local systems, occur in $\pi_!(1)$. But then, since each complex occuring lies in the canonical basis, Lusztig's description of this set implies that if $j: E_V^r \to E_V$ then $j^*(\pi_!(1))$ must be precisely
\[
\bigoplus_{\mu, \lambda} K_{\mu, \lambda} j^*(IC(\mathcal L_\lambda)),
\]
(this follows more directly from the fact that the map $\pi$ is small over the regular set). Now \cite[Theorem 5.2]{L98} shows that the $b_\lambda$ are obtained from the trace of Frobenius on the stalks of $IC(\mathcal L_\lambda)$. Finally, let $\tilde Y$ be the preimage of the set of regular elements, and let $\tilde{\pi}$ be the restriction of $\pi$ to $\tilde Y$. Then by base change we have $j^*(\pi_!(1)) = \tilde{\pi}_!(1_{\tilde Y})$.  An element of $(x, \mathfrak f) \in Y$ lies in $\tilde Y$ precisely when all the representations induced by $x$ on the associated graded of the filtration are regular (this is of course the reason the subspace $\mathcal F^r$ is a subalgebra as mentioned before). It is easy to pave the fibers of $\tilde \pi$ by affine spaces, and so it follows that they are pure, hence we see that the trace of Frobenius on $\tilde{\pi}_!(1_{\tilde Y})$ is $\rho_\lambda$. The result follows.
\end{proof}

We now easily obtain the description of the $s_\lambda$ in $\mathcal F \cong \mathbf U^+_v$:

\begin{cor}
\label{imaginary}
For all partitions $\lambda$, we have
\[
s_\lambda = r(b_\lambda).
\]
Thus $\Uz$ at $v = \e$ corresponds under $\chi_\e$ to $\mathcal F^r$, and $\mathcal F^r$
is precisely $r(\mathcal F)$.
\end{cor}
\begin{proof}
 Using the above identification of $\Uz$ with the algebra of symmetric functions
we see (say from \cite[\S 1.6]{M}) that
\begin{equation}
\tilde{P}_\lambda = \bigoplus_{\mu \geq \lambda} K_{\mu \lambda} s_{\mu}.
\end{equation}
where the $K_{\mu \lambda}$ are the Kostka numbers, exactly the relation between the
$\tilde{P}_\lambda$ and $j^*(b_\lambda)$ in the last proposition. Since the matrix of
Kostka polynomials is invertible we immediately get the result. The last sentence follows
immediately from the description of the canonical basis.
\end{proof}

Finally, we can restate this corollary in purely algebraic fashion. To do this we need some of the set up of the \cite{BCP}. In that paper they define subalgebras of $\UA$, denoted $\mathbf U^+(>)$ and $\mathbf U^+(<)$ (in fact they use the notation
$\mathbf U^+(>)_{\A}$ etc. but we need only the integral form here, so we abbreviate for
simplicity of notation), which are the subalgebras generated
by $(E_{1,k})_{k \geq 1}$ and $(E_{0,k})_{k \geq 1}$ respectively.
Then they show that there is a direct sum decomposition
\begin{equation}
\label{direct}
\UA = \mathbf U^+(0) \oplus \biggl( \mathbf U^+(>)\mathbf U^+(0)\mathbf U^+(<)_{+} +
\mathbf U^+(>)_{+}\mathbf U^+(0)\mathbf U^+(<)\biggr).
\end{equation}

Let $\pi^0$ denote the projection
onto $\mathbf U^+(0)$ associated to this decomposition. Then we have

\begin{cor}
For all partitions $\lambda$ we have
\[
\pi^0(b_\lambda) = s_{\lambda}.
\]
\end{cor}
\begin{proof}
From Remark \ref{generic} we know that the generic situation can be recovered from an infinite family of specializations. Hence the result is clear once we identify $\pi^0$ with the map $r: \mathcal F \to \mathcal F^r$. It is known from \cite{BCP} that the direct sum in \eqref{direct} is orthogonal. Let $(\cdot,\cdot)_V$ denote the inner product on $\mathcal F_V$ given by
\[
(f,g) = q^{dim(G_V)}|G_V^F|^{-1} \sum_{x \in E_V^F} f(x)g(x).
\]
Taking the orthogonal direct sum we obtain an inner product on $\mathcal F$, which is
shown in \cite{L98} to coincide under $\chi_v$ with the inner product on $\UA$. This shows
immediately that orthogonal projection onto $\mathbf U^+(0)$ is precisely the map $r$, and
the result follows.
\end{proof}

\section{The Drinfeld Presentation}

We wish to observe in this section that it is now easy to check, for the plus part of
the quantum group, that the Drinfeld presentation of $\mathbf U^+$ can be realized in the
context of the Hall algebra. Of course, in a sense all this represents is a translation of the insight of Beck in \cite{B} into the quiver context, or alternatively a marriage of the work of Kac \cite{Kac} and Ringel \cite{R}. Recall from section \ref{quantum} that $\mathbf U^+$ is generated by the elements $\{x_r^+, C^{s}k^{-1}x_s^-, h_s \colon r\geq 0, s > 0\}$. For simplicity of notation, we will write $y_s^- = C^{s}k^{-1}x^-_s$. Then we have the relations:

\begin{enumerate}
\item $[h_s,x_r^+] = \pm \frac{[2s]}{s} x_{r+s}^{+}$
\item $[h_s,y_r^-] = \pm \frac{[2s]}{s} y_{r+s}^-$
\item $x_{r+1}^+x_s^+ - v^{2}x_s^+x_{r+1}^+ = v^{2}x_r^+x_{s+1}^+ -
x_{s+1}^+x_r^+$,
\item $y_{r+1}^-y_s^- - v^{-2}y_s^-y_{r+1}^- = v^{-2}y_r^-y_{s+1}^- -
y_{s+1}^-y_r^-$,
\item $v^{-2}x_r^+y_s^- - y_s^-x_r^+ = \psi_{r+s}^+/(v-v^{-1})$.
\end{enumerate}
where $\psi_i = C^{i/2}k^{-1}\psi^+_i$. (Of course, we are here assuming that this algebra injects into the algebra defined in section \ref{quantum}, however this is known to be true, see \cite{B}). Let $\mathbf U^+_D$ denote this algebra.

We define elements $\eta_r$ of $\mathcal F$ by the generating function equation
\[
\sum_{r \geq 0} \rho_r u^r = \rm{exp}\biggl(\sum_{r\geq 1} \frac{\eta_r u^r}{[r]_\e}\biggr).
\]
Explicitly we have
\[
\rho_k = \frac{1}{k}\sum_{s=1}^{k} \frac{s}{[s]_\e}\eta_s\rho_{k-s}.
\]
Let $\mathcal S$ be as in Remark \ref{generic}, and note that we have natural elements $\eta_s, \rho_k, \mu_r, \gamma_r$ in $\mathcal F_\mathcal S$.

\begin{prop}
There is an isomorphism $\Lambda: \mathbf U_D^+ \to \mathcal F_\mathcal S$ given by $x_r^+ \mapsto \mu_{r+1}$, $y_r^- \mapsto -\gamma_r$ and $h_r \mapsto \eta_r$.
\end{prop}
\begin{proof}
To show that $\Lambda$ is well defined we must check the relations $(1)$ to $(5)$. It is enough to check in each specialization separately. Now relation $(5)$ is just $(3)$ in Corollary \ref{properties}, and so has already been checked. Relations $(3)$ and $(4)$ are straightforward to check, so the only difficulty is to show $(1)$ and $(2)$. An argument similar to the proof of Theorem \ref{relation} shows the following.

\begin{lemma}
\label{comm}
For $r,s \geq 1$ we have
\begin{enumerate}
\item $\rho_r\mu_s = \sum_{i=0}^r [r-i+1]_\e \mu_{r+s-i}\rho_i$,
\item $\gamma_s\rho_r = \sum_{i=0}^r [r-i+1]_\e \rho_i\gamma_{r+s-i}$.
\end{enumerate}
\end{lemma}
\begin{proof}
By using the action of $T_\varpi$ we may assume that $s=1$. Then $(1)$ follows by observing that a module $M$ which contains a submodule $N$ isomorphic to $S_1$ with $M/N$ regular of dimension $(r,r)$ must by a dimension count and \eqref{homs} be isomorphic to a $R \oplus P_k$ where $k \leq r$ and $R$ is regular. One can prove $(2)$ by a similar argument, or appeal to duality.
\end{proof}

It is straightforward (though at least for the author, somewhat painful), to check that relations $(1)$ and $(2)$ are equivalent to those in Lemma \ref{comm} using induction and the identity
\[
\sum_{i=1}^m \frac{[2i]}{[i]}[m-i+1] = m[m+1], \qquad m\geq 1.
\]

Thus $\Lambda$ is a well defined homomorphism. That it is an isomorphism in the generic situation (see \cite{L98} for more details about this) can be proved by exhibiting a basis
of both algebras.
\end{proof}

\begin{remark}
The point here is just to show that Drinfeld-style generators for the algebra arise naturally in this context -- though for us the Heisenberg generators are not transparent as they do not lie in the integral form, which is perhaps also why it is inconvenient to prove those relations. (The relations in Lemma \ref{comm} are proved in \cite{CP} starting from the relations of the Drinfeld presentation.) For us the integral generators $x^+_r, y^-_s, \tilde{P}_k$ are the easiest to describe. Since the whole quantum group can be constructed, via the Drinfeld double, from $\mathbf U^+$ and its inner product, one could presumably recover the whole presentation in this way, once certain inner products and coproducts were known.
\end{remark}

We wish to make one final observation. The elements $\tilde{P}_k$ determine the Drinfeld polynomials attached to a finite dimensional irreducible representation of the affine quantum group. It is important in applications to understand the behaviour of these elements under the coproduct. Let $\Delta$ be the twisted coproduct described in Section \ref{Kronecker}, or rather its analogue for $\mathcal F$. Define $\mathcal F^>, \mathcal F^<$ and $\mathcal F^r$ to be the subalgebras of $\mathcal F$ generated by the $(\mu_k), (\gamma_k)$, and $(\rho_k)$ respectively. Under $\chi_v$ these correspond to the subalgebras $\mathbf U^+(>), \mathbf U^+(<)$ and $\mathbf U^+(0)$ of \cite{BCP}. It is not hard to see that $\mathcal F = \mathcal F^r \oplus \mathcal F^>_+\mathcal F \oplus \mathcal F\mathcal F^<_+$. Also for any algebra $R$ let $R_+$ denote the augmentation algebra.

\begin{lemma}
We have for $k \geq 0$
\[
\Delta(\rho_k) = \sum_{i=0}^k \rho_i \otimes \rho_{k-i} + \text{functions in } \mathcal F^<_+ \mathcal F^r \otimes \mathcal F^r \mathcal F^>_+
\]
and hence we have 
\[
\Delta(\tilde{P}_k) = \sum_{i=0}^k \tilde{P}_i \otimes \tilde{P}_{k-i} + \mathbf U^+(<)_+ \mathbf U^+(0) \otimes \mathbf U^+(0) \mathbf U^+(>)_+.
\]
\end{lemma}
\begin{proof}
There are two points to note. The first is to observe the following: if $W < V$ are vector spaces of dimensions $(k,k)$ and $(n,n)$ respectively, and $(V,x)$ is a representation of $K$ such that $x$ preserves $W$ then $(V,x)$ is regular if and only if $(W,x)$ and $(V/W, \overline x)$ are both regular. In order to see that the other terms are as claimed one must note that a regular module cannot contain a preinjective submodule or have a preprojective extension. The formula in $\UA$ follows by passing to the generic situation as in Remark \ref{generic}.
\end{proof}

\section{Generalizations}

The elements $b_\lambda$ can be thought of as the purely
imaginary part of the canonical basis. It would be desirable to have a generalization
of the results this paper to arbitrary affine quantum groups, or at least the symmetric
simply laced ones (note that the results of \cite{BCP} have been generalized to arbitrary affine quantum groups in \cite{BN}). The key seems to be to understand in geometric terms the meaning of an imaginary root (I thank George Lusztig for clarifying this for me).
For $\widehat{\mathfrak{sl}}_2$ this is easy, as a root is imaginary precisely 
when it is regular (see \cite{L92} for this terminology). However in other types
this is not the case -- an imaginary root is always regular, but in general there
are real regular roots. 

In the algebraic context of \cite{B} one passes from $\asl$ to other types by means
of ``vertex'' embeddings of $\asl$ -- one for each node of the finite type diagram.
The meaning of these embeddings in the context of quiver representations seems straightforward for a node of valence one, however, say for the central node
in $\tilde{D}_4$, it is not clear to the author what one should do. However,
for all but this central node, there seems to be a natural reduction that we
will now describe.

We work in the context of the McKay correspondence as in \cite{L92}.
Let $\Gamma$ be a finite subgroup of $\mathrm{SL}_2(\mathbb C)$, which we will
assume contains $-I$, and consider the action of $\Gamma$ on $\mathbb P^1$. 
There is a natural notion of a regular representation of the associated quiver,
and such representations have a spectrum which consists of orbits of the action of $\Gamma$
on $\mathbb P^1$ (for $\asl$ we have $\Gamma = \{\pm I\}$, and so the orbits 
are all single points). The points $L$ of $\mathbb P^1$ which have stabilizer
$\Gamma_L$ of order greater than two are where the subtlety lies. If $\Xi$ denotes
a set of representatives for the orbits of $\Gamma$ on these points, then there is
a correspondence between $\Xi$ and the finite type Dynkin diagram as follows.
Each element of $\Xi$ corresponds to an ``arm'' of the Dynkin diagram (so $\Xi$
has two or three elements), and the length of the arm associated to $L \in \Xi$
is $1/2|\Gamma_L|$ (including the central node). The representations which have
spectrum in the $\Gamma$-orbit of $L$ are equivalent to the category of nilpotent representations of a cyclic quiver of order $1/2|\Gamma_L|$. Thus using these 
embedding we can reduced to finding vertex embeddings for $\widehat{\mathfrak{sl}}_n$
and, separately, for the central node in types $D$ and $E$.

For the central node, it seems reasonable to attach the canonical basis
elements whose supports contain generic regular representations (i.e. those
whose spectrum consists of points in $\mathbb P^1$ with minimal stabilizer).
In Lusztig's description of the canonical basis, these are the ones which 
come from symmetric group covering, so already one sees symmetric functions.


\textit{Acknowledgements.} 
I would like to thank George Lusztig for useful conversations, and Anthony Henderson for
a careful critique of a draft of the paper.

\end{document}